# STATISTICAL ROMBERG EXTRAPOLATION: A NEW VARIANCE REDUCTION METHOD AND APPLICATIONS TO OPTION PRICING

By Ahmed Kebaier

*Université de Marne-La-Vallée*

We study the approximation of $\mathbb{E}f(X_T)$ by a Monte Carlo algorithm, where $X$ is the solution of a stochastic differential equation and $f$ is a given function. We introduce a new variance reduction method, which can be viewed as a statistical analogue of Romberg extrapolation method. Namely, we use two Euler schemes with steps $\delta$ and $\delta^\beta, 0 < \beta < 1$. This leads to an algorithm which, for a given level of the statistical error, has a complexity significantly lower than the complexity of the standard Monte Carlo method. We analyze the asymptotic error of this algorithm in the context of general (possibly degenerate) diffusions. In order to find the optimal $\beta$ (which turns out to be $\beta = 1/2$), we establish a central limit type theorem, based on a result of Jacod and Protter for the asymptotic distribution of the error in the Euler scheme. We test our method on various examples. In particular, we adapt it to Asian options. In this setting, we have a CLT and, as a by-product, an explicit expansion of the discretization error.

**1. Introduction.** In many numerical problems—in particular, in mathematical finance—one has to compute, using the Monte Carlo method, $\mathbb{E}f(X_T)$, where $X_T$ is a diffusion process. The advantage of using a probabilistic approach instead of a PDE approach is that one may solve problems in high dimension. But on the other hand, the Monte Carlo algorithms are much slower and their practical efficiency highly depends on the variance of the random variables at hand. This is why variance reduction plays a crucial role in the practical implementation.

There are several classes of methods which are used to reduce variance: the control variate approach, the antithetic variate method, moment matching, importance sampling, conditional Monte Carlo methods.... (For more details about variance reduction methods, see [3].) In this paper we intro-
duce a new variance reduction method which we will call statistical Romberg
method. This method can be viewed as a control variate method. Roughly
speaking, the idea is as follows: use many sample paths with a coarse time
discretization step and few additional sample paths with a fine time dis-
cretization step. More precisely, in order to construct our control variate, we
discretize the diffusion $(X_t)_{0 \le t \le T}$ by two Euler schemes with time steps $T/n$
and $T/m$ ($m \ll n$). Suppose for a while that $M_m := \mathbb{E}f(X_T^m)$ is known. Then
we construct

$$Q = f(X_T^n) - f(X_T^m) + M_m.$$

Clearly, we have $\mathbb{E}f(X_T^n) = \mathbb{E}Q$ and we wish to simulate $Q$ instead of $f(X_T^n)$.
Standard arguments show that for $f$ Lipschitz continuous, we have

$$\mathrm{Var}(Q) = O(1/m)$$

(see Proposition 3.1). So the variance of $Q$ tends to zero as $m$ tends to infinity
and, consequently, in order to achieve a given accuracy, we need a much
smaller sample. This significantly reduces the complexity of the algorithm.
But, on the other hand, we have to compute the quantity $M_m$, and this is
done again by Monte Carlo sampling. This will increase the complexity of
the algorithm and we need to find a good balance which guarantees that
the global complexity decreases. Let $N_m$ be the number of Monte Carlo
simulations used for the evaluation of $M_m$ and $N_n$ be the number of Monte
Carlo simulations used for the evaluation of $\mathbb{E}Q$. The question is how to
choose $m$, $N_m$ and $N_n$.

In the classical Monte Carlo method, one needs to choose the number $n$
of intervals in the discretization and the number $N$ of simulations. The
parameter $n$ drives the so-called *discretization error* due to discretization,
whereas the number $N$ controls the *statistical error*. For a rational choice
of $N$ versus $n$, one may try to minimize the error for a given computational
time (see [5]) or, equivalently, to minimize the computational time for a
given total error (see [12]).

An optimal choice of the number of Monte Carlo samples must be based
on a precise evaluation of the discretization error. This requires some kind
of regularity: in [2], the regularizing effect of the diffusion process allows
to prove that the rate of convergence is $1/n$ for measurable functions. A
Hörmander type assumption is needed. In the case of general diffusions, the
same order of convergence can be proved for regular functions only. (See the
pioneering paper [16] for $\mathscr{C}^6$-functions and for $\mathscr{C}^3$-functions, see [12]. For
stochastic differential equations driven by a Lévy process, see [8].)

Our first result shows that, in the context of possibly degenerate dif-
fusions, the discretization error for $\mathscr{C}^1$-functions is at least $o(1/\sqrt{n})$ (see

Proposition 2.2). For such functions, we give an example for which the discretization error is of order $1/n^\alpha$ for any $\alpha \in (1/2, 1]$ (see Proposition 2.3).

Our second result is a central limit theorem for the statistical Romberg method (see Theorem 3.2). This theorem uses the weak convergence of the normalized error of the Euler scheme for diffusions proved by Kurtz and Protter [11] (and strengthened by Jacod and Protter [9]).

Based on this central limit theorem, we are able to fix the optimal balance between $m$, $N_m$ and $N_n$. It turns out that, for a given error level $\varepsilon = 1/n^\alpha$, we obtain $m = \sqrt{n}$, $N_m = n^{2\alpha}$ and $N_n = n^{2\alpha-1/2}$. With this choice, the complexity of our algorithm is $C_{\mathrm{SR}} = C \times n^{2\alpha+1/2}$, $C > 0$, while the complexity in the standard Monte Carlo is $C_{\mathrm{MC}} = C \times n^{2\alpha+1}$. So we have a clear gain.

Our approach can also be used for the Monte Carlo approximation of expectations of functionals of the whole path of a diffusion. In particular, we investigate the application to the pricing of Asian options. In this setting, the approximation relies on the discretization of the integral of the price process over a time interval. It was shown in [17] that the trapezoidal rule is one of the most efficient methods for this discretization. We analyze the error process, which is of the order $1/n$. We prove a stable functional central limit theorem in the spirit of Jacod and Protter [9]. As a consequence of this result, we give an expansion of the analytical error, which, in contrast with Temam's [17] result, does not use the associated PDE We also use our result in order to optimize the choice of the parameters, which are different from the ones in the Euler scheme.

The organization of the paper is the following. We first recall essential facts about the Euler scheme, including error evaluations. In Section 3 we describe our statistical Romberg method for the Euler scheme. In Section 4 we apply the idea of statistical Romberg approximation to the discretization of path integrals of a diffusion in the context of Asian option pricing. Section 5 is devoted to numerical tests and comparisons. In the last section we give conclusions and remarks.

## 2. On the discretization error of the Euler scheme.

Let $(X_t)_{0 \le t \le T}$ be the process with values in $\mathbb{R}^d$, solution to

$$(1) \qquad dX_t = b(X_t)\, dt + \sigma(X_t)\, dW_t, \qquad X_0 = x \in \mathbb{R}^d,$$

where $W = (W^1, \ldots, W^q)$ is a $q$-dimensional Brownian motion on some given filtered probability space $\mathcal{B} = (\Omega, \mathcal{F}, (\mathcal{F}_t)_{t \ge 0}, P)$. $(\mathcal{F}_t)_{t \ge 0}$ is the standard Brownian filtration. The functions $b : \mathbb{R}^d \longrightarrow \mathbb{R}^d$ and $\sigma : \mathbb{R}^d \longrightarrow \mathbb{R}^{d \times q}$ are continuously differentiable and satisfy $\exists C_T > 0; \forall x, y \in \mathbb{R}^d$, we have

$$|b(x) - b(y)| + |\sigma(x) - \sigma(y)| \le C_T |y - x|.$$

We consider the Euler continuous approximation $X^n$ with step $\delta = T/n$ given by

$$dX_t^n = b(X_{\eta_n(t)})\, dt + \sigma(X_{\eta_n(t)})\, dW_t, \qquad \eta_n(t) = [t/\delta]\delta.$$

It is well known that the Euler scheme satisfies the following properties (see, e.g., [6]):

$$(2) \qquad \forall p > 1 \qquad \mathbb{E} \sup_{t \in [0,T]} |X_t - X_t^n|^p \le \frac{K_p(T)}{n^{p/2}}, \qquad K_p(T) > 0.$$

$$(3) \qquad \forall p > 1 \qquad \mathbb{E} \sup_{t \in [0,T]} |X_t|^p + \mathbb{E} \sup_{t \in [0,T]} |X_t^n|^p \le K_p'(T), \qquad K_p'(T) > 0.$$

2.1. *Stable convergence and the Euler scheme error.* We first recall basic facts about stable convergence. In the following we adopt the notation of Jacod and Protter [9]. Let $X_n$ be a sequence of random variables with values in a Polish space $E$, all defined on the same probability space $(\Omega, \mathcal{F}, P)$. Let $(\tilde{\Omega}, \tilde{\mathcal{F}}, \tilde{P})$ be an extension of $(\Omega, \mathcal{F}, P)$, and let $X$ be an $E$-valued random variable on the extension. We say that $(X_n)$ converges in law to $X$ stably and write $X_n \overset{\text{stably}}{\Longrightarrow} X$, if

$$\mathbb{E}\left(Uh(X_n)\right) \to \tilde{\mathbb{E}}\left(Uh(X)\right),$$

for all $h : E \to \mathbb{R}$ bounded continuous and all bounded random variable $U$ on $(\Omega, \mathcal{F})$. This convergence, introduced by Rényi [15] and studied by Aldous and Eaglson [1], is obviously stronger than convergence in law. The following lemma will be crucial:

If $V$ is another variable with values in another Polish space $F$, we have the following result:

LEMMA 2.1. *Let $V_n$ and $V$ be defined on $(\Omega, \mathcal{F})$ with values in another metric space $E'$.*

$$\text{If } V_n \overset{\mathbb{P}}{\to} V, X_n \overset{\text{stably}}{\Longrightarrow} X, \qquad \text{then } (V_n, X_n) \overset{\text{stably}}{\Longrightarrow} (V, X).$$

*This result remains valid when $V_n = V$.*

For a proof of this lemma, see [9].

Note that all this applies when $X_n$, $X$ are $\mathbb{R}^d$-valued càdlàg processes, with $E = \mathbb{D}([0,T], \mathbb{R}^d)$. Now assume that

$$\varphi(X_t) = \begin{pmatrix} b_1(X_t) & \sigma_{11}(X_t) & \dots & \sigma_{1q}(X_t) \\ b_2(X_t) & \sigma_{21}(X_t) & \dots & \sigma_{2q}(X_t) \\ \vdots & \vdots & & \vdots \\ b_d(X_t) & \sigma_{d1}(X_t) & \dots & \sigma_{dq}(X_t) \end{pmatrix} \quad \text{and} \quad dY_t := \begin{pmatrix} dt \\ dW_t^1 \\ \vdots \\ dW_t^q \end{pmatrix},$$

then the SDE (1) becomes

$$(2) \qquad\qquad\qquad dX_t = \varphi(X_t) \, dY_t.$$

The Euler continuous approximation $X^n$ with step $\delta = T/n$ is given by

$$dX^n_t = \varphi(X_{\eta_n(t)})\,dY_t, \qquad \eta_n(t) = [t/\delta]\delta.$$

The following result proven by Jacod and Protter [9] is an improvement on the result given by Kurtz and Protter [11].

THEOREM 2.1. *With the above notation we have*

$$\sqrt{n}U^n =: \sqrt{n}(X^n - X) \overset{stably}{\Longrightarrow} U,$$

*with $U$ a $d$-dimensional process satisfying*

$$(3) \qquad dU^i_t = \sum_{j=1}^{q+1}\sum_{k=1}^{d} \varphi'^{ij}_k(X_t)\left[U^k_t\,dY^j_t - \sum_{l=1}^{q+1}\varphi^{kl}(X_t)\,dN^{lj}_t\right], \qquad U^i_0 = 0$$

*($\varphi'^{ij}_k$ is the partial derivative of $\varphi^{ij}$ with respect to the $k$th coordinate), and $N$ is given by*

$$N^{1i} = 0, \qquad 1 \le i \le q+1,$$
$$N^{j1} = 0, \qquad 1 \le j \le q+1,$$
$$N^{ij} = \frac{B_{ij}}{\sqrt{2}}, \qquad 2 \le i,j \le q+1,$$

*where $(B^{ij})_{1 \le i,j \le q}$ is a standard $(q)^2$-dimensional Brownian motion defined on an extension $\tilde{\mathcal{B}} = (\tilde{\Omega}, \tilde{\mathcal{F}}, (\tilde{\mathcal{F}}_t)_{t \ge 0}, \tilde{P})$ of the space $(\Omega, \mathcal{F}, (\mathcal{F}_t)_{t \ge 0}, P)$, which is independent of $W$.*

We will need the following property of the process $U$.

PROPOSITION 2.1. *Under the assumptions of the above theorem, we have*

$$\tilde{\mathbb{E}}(U_T/\mathcal{F}_T) = 0.$$

PROOF. Consider the unique solution of the $d$-dimensional linear equation

$$(4) \qquad \mathcal{E}_T = I_d + \sum_{j=1}^{q+1} \int_0^T \varphi'^j(X_t)\mathcal{E}_t\,dY^j_t,$$

where $\varphi'^j$ is a $d \times d$ matrix with $(\varphi'^j)_{ik} = \varphi'^{ij}_k$. From Theorem 56, page 271, in [14], it follows that

$$(5) \qquad U_T = -\sum_{j=1}^{q+1} \mathcal{E}_T \int_0^T \mathcal{E}_t^{-1}\varphi'^j(X_t)\varphi(X_t)\,dN^j_t.$$

If $Z$ is a bounded $\mathcal{F}_T$-measurable r.v., we have

$$\mathbb{E}\,(U_T \cdot Z) = -\sum_{j=1}^{q+1} \mathbb{E}\left( Z \mathcal{E}_T \int_0^T \mathcal{E}_t^{-1} \varphi'^j(X_t)\varphi(X_t)\,dN_t^j\right) = 0,$$

as can be seen by representing $Z\mathcal{E}_T$ as a stochastic integral w.r.t. $W$. $\quad\square$

2.2. *The discretization error.*  In the following we focus on the discretization error given by the bias

(6) $$\varepsilon_n := \mathbb{E}f(X_T^n) - \mathbb{E}f(X_T),$$

where $f$ is a given function. Talay and Tubaro [16] prove that if $f$ is sufficiently smooth, then $\varepsilon_n \sim c/n$ with $c$ a given constant. A similar result was proven by Kurtz and Protter [12] for a function $f \in \mathscr{C}^3$. The same result was extended in [2] for a measurable function $f$, but with a nondegeneracy condition of Hörmander type on the diffusion. In the context of possibly degenerate diffusions, the discretization error for functions which are not $\mathscr{C}^3$ is not yet completely understood. For a Lipschitz-continuous function $f$, the estimate $|\varepsilon_n| \leq \frac{c}{\sqrt{n}}$ follows easily from  (2). The following proposition and the example below focus on the rate of convergence of $\varepsilon_n$ for $\mathscr{C}^1$ functions.

PROPOSITION 2.2.  *Let $f$ be an $\mathbb{R}^d$-valued function satisfying*

(7) $$|f(x) - f(y)| \leq C(1 + |x|^p + |y|^p)|x - y| \qquad \text{for some } C, p > 0.$$

*Assume that $\mathbb{P}(X_T \notin \mathcal{D}_f) = 0$, where $\mathcal{D}_f := \{x \in \mathbb{R}^d | f$ is differentiable at $x\}$, then*

$$\lim_{n \to \infty} \sqrt{n}\,\varepsilon_n = 0.$$

PROOF.  We have, with probability 1,

$$\sqrt{n}(f(X_T^n) - f(X_T)) = \sqrt{n}\,\nabla f(X_T) \cdot U_T^n + R_n,$$

with

$$R_n = \sqrt{n}|U_T^n|\varepsilon(X_T, U_T^n) \quad \text{and} \quad \varepsilon(X_T, U_T^n) \xrightarrow{\mathbb{P}} 0.$$

It follows that $R_n \xrightarrow{\mathbb{P}} 0$, since $(\sqrt{n}|U_T^n|)$ is tight. Consequently, we deduce, using (13), (2), (3) and Theorem 2.1, that

$$\sqrt{n}\,\varepsilon_n \to \mathbb{E}(\nabla f(X_T) \cdot U_T),$$

and using Proposition 2.1, it follows that

$$\mathbb{E}(\nabla f(X_T) \cdot U_T) = 0,$$

which completes the proof. $\quad\square$

The following example proves that, for $\alpha \in (1/2, 1]$, there exists a $\mathscr{C}^1$ function with bounded derivatives and a diffusion $X$ such that

$$(7) \qquad n^\alpha \varepsilon_n \to C_f(T, \alpha),$$

where $C_f(T, \alpha)$ is positive. In other words, the rate of convergence can be $1/n^\alpha$ for all values of $\alpha \in (1/2, 1]$.

EXAMPLE 2.1.   Consider the bi-dimensional diffusion $Z = (X, Y)$ satisfying the following SDE:

$$(8) \qquad \begin{aligned} dX_t &= -X_t/2 \, dt - Y_t \, dW_t, \\ dY_t &= -Y_t/2 \, dt + X_t \, dW_t \end{aligned}$$

and the map $f_\alpha : z = (x, y) \mapsto ||z||^2 - 1|^{2\alpha}$. The solution of (8), subject to $Z_0 = (\cos\theta, \sin\theta)$, is given by $Z_t = (\cos(\theta + W_t), \sin(\theta + W_t))$. We assume that $\theta \in [0, 2\pi]$, so the diffusion $Z$ lives on the unit circle.

PROPOSITION 2.3.   *Let $Z^n$ be the Euler scheme associated with $Z$. For $\alpha \in [1/2, 1]$, we have*

$$(9) \qquad n^\alpha \mathbb{E}(f_\alpha(Z_t^n) - f_\alpha(Z_t)) \to (2t)^\alpha \mathbb{E}|G|^{2\alpha}, \qquad t \geq 0,$$

*where $G$ is a standard normal r.v.*

PROOF.   We have, since $f_\alpha$ vanishes on the unit circle,

$$\begin{aligned} n^\alpha \mathbb{E}(f_\alpha(Z_t^n) - f_\alpha(Z_t)) &= n^\alpha \mathbb{E}||Z_t^n|^2 - 1|^{2\alpha} \\ &= n^\alpha \mathbb{E}|[Z_t + Z_t^n] \cdot [Z_t^n - Z_t]|^{2\alpha}. \end{aligned}$$

Using Theorem 2.1,

$$(10) \qquad n^\alpha \mathbb{E}(f_\alpha(Z_t^n) - f_\alpha(Z_t)) \to 2^{2\alpha} \mathbb{E}|Z_t \cdot \tilde{U}_t|^{2\alpha},$$

where $\tilde{U} = (\tilde{U}^1, \tilde{U}^2)$ is given by

$$(11) \qquad \begin{aligned} d\tilde{U}_t^1 &= -\tfrac{1}{2}\tilde{U}_t^1 \, dt - \tilde{U}_t^2 \, dW_t + \tfrac{1}{\sqrt{2}} X_t \, d\tilde{B}_t \\ d\tilde{U}_t^2 &= -\tfrac{1}{2}\tilde{U}_t^2 \, dt + \tilde{U}_t^1 \, dW_t + \tfrac{1}{\sqrt{2}} Y_t \, d\tilde{B}_t \end{aligned}$$

and $\tilde{B}$ is a standard Brownian motion independent of $W$. The solution of (11) is given by

$$(12) \qquad \begin{aligned} \tilde{U}_t^1 &= \tfrac{1}{\sqrt{2}} X_t \tilde{B}_t \\ \tilde{U}_t^2 &= \tfrac{1}{\sqrt{2}} Y_t \tilde{B}_t, \end{aligned}$$

which completes the proof.   $\square$

**3. The statistical Romberg method.** Before introducing our algorithm, we recall some essential facts about the Monte Carlo method. In many applications (in particular, for the pricing of financial securities), the effective computation of $\mathbb{E}\, f(X_T)$ is crucial (see, e.g. [13]). The Monte Carlo method consists of the following steps:

- Approximate the process $(X_t)_{0 \le t \le T}$ by the Euler scheme $(X_t^n)_{0 \le t \le T}$, with step $T/n$, which can be simulated.
- Evaluate the expectation on the approximating process $f(X_T^n)$ by the Monte Carlo method.

In order to evaluate $\mathbb{E} f(X_T^n)$ by the Monte Carlo method, $N$ independent copies $f(X_{T,i}^n)_{1 \le i \le N}$ of $f(X_T^n)$ are sampled and the expectation is approximated by the following quantity:

$$\hat{f}^{n,N} := \frac{1}{N} \sum_{i=1}^{N} f(X_{T,i}^n).$$

The approximation is affected by two types of errors. The discretization error $\varepsilon_n$, studied in the above section, and the statistical error $\hat{f}^{n,N} - \mathbb{E} f(X_T)$, controlled by the central limit theorem and which is of order $1/\sqrt{N}$. An interesting problem (studied in [5] and [12]) is to find $N$ as a function of $n$ so that both errors are of the same order.

The following result highlights the behavior of the global error in the classical Monte Carlo method. It can be proved in the same way as the limit theorem given in [5].

THEOREM 3.1. *Let $f$ be an $\mathbb{R}^d$-valued function satisfying*

$$(13) \qquad |f(x) - f(y)| \le C(1 + |x|^p + |y|^p)|x - y| \qquad \text{for some } C, p > 0.$$

*Assume that $\mathbb{P}(X_T \notin \mathcal{D}_f) = 0$, where $\mathcal{D}_f := \{x \in \mathbb{R}^d | f \text{ is differentiable at } x\}$, and that for some $\alpha \in [1/2, 1]$, we have*

$$(14) \qquad \lim_{n \to \infty} n^\alpha \varepsilon_n = C_f(T, \alpha).$$

*Then*

$$n^\alpha \left( \frac{1}{n^{2\alpha}} \sum_{i=1}^{n^{2\alpha}} f(X_{T,i}^n) - \mathbb{E}\, f(X_T) \right) \Longrightarrow \sigma \bar{G} + C_f(T, \alpha),$$

*with $\sigma^2 = \text{Var}(f(X_T))$ and $\bar{G}$ a standard normal.*

A functional version of this theorem, with $\alpha = 1$, was proven by Kurtz and Protter [12] for a function $f$ of class $\mathscr{C}^3$. We can interpret the theorem as follows. For a total error of order $1/n^\alpha$, the minimal computation effort necessary to run the Monte Carlo algorithm is obtained for $N = n^{2\alpha}$. This leads to an optimal time complexity of the algorithm given by

$$(13) \quad C_{\mathrm{MC}} = C \times (nN) = C \times n^{2\alpha+1} \qquad \text{with } C \text{ some positive constant.}$$

Recall that the time complexity of an algorithm $A$ is proportional to the maximum number of basic computations performed by $A$.

3.1. *The Euler scheme and the statistical Romberg method.* It is well known that the rate of convergence in the Monte Carlo method depends on the variance of $f(X_T^n)$, where $X_T^n$ is the Euler scheme of step $T/n$. This is a crucial point in the practical implementation. A large number of reduction of variance methods are used in practice. Our algorithm proposes a control variate reduction of variance. Its specificity is that the control variate is constructed itself using the Monte Carlo method, applied to the same discretization scheme, but with a step $m$ which is specifically lower than the approximation step $n$ (using two discretization steps is an idea which already appears in Romberg's method). Let us be more precise. We fix $m \ll n$ and we denote

$$Q = f(X_T^n) - f(X_T^m) + M_m,$$

where $M_m = \mathbb{E}f(X_T^m)$ and we suppose for a while that $M_m$ is known. Note that

$$\mathbb{E}(Q) = \mathbb{E}f(X_T^n),$$

so that $f(X_T^m) - M_m$ appears as a control variate.

Consider a function $f : \mathbb{R}^d \longrightarrow \mathbb{R}^d$ which is Lipschitz continuous of constant $[f]_{\mathrm{lip}}$, that is, $[f]_{\mathrm{lip}} = \sup_{u \neq v} \frac{|f(u) - f(v)|}{|u - v|}$.

PROPOSITION 3.1. *Under the above assumptions, we have*

$$(14) \quad \sigma_Q^2 := \mathrm{Var}(Q) = O(1/m).$$

PROOF. We have

$$\begin{aligned}
\sigma_Q &= \|Q - \mathbb{E}Q\|_2 \\
&\leq \|f(X_T^n) - f(X_T^m)\|_2 \\
&\leq [f]_{\mathrm{lip}} \Big[ \sup_{t \in [0,T]} \|X_t - X_t^n\|_2 + \sup_{t \in [0,T]} \|X_t - X_t^m\|_2 \Big].
\end{aligned}$$

Using (2), we deduce that $\exists K' > 0$ such that

$$\sigma_Q \leq K' \left( \frac{1}{\sqrt{m}} \right),$$

which completes the proof.  $\square$

Inequality (14) shows that the variance of $Q$ is significantly smaller than the variance of $f(X_T^n)$, so that $Q$ appears as a good candidate for the reduction of variance method. However, computing $\mathbb{E}Q$ supposes to compute $M_m = \mathbb{E}f(X_T^m)$ in the first place, and this is also done by the Monte Carlo method. So there is a certain extra quantity of computation to be done. In practice, the sample paths used for the computation of $M_m$ will be independent of those used for the computation of $Q$.

In the following we make a complexity analysis which permits to choose $m$ as a function of $n$ in order to minimize the complexity of the algorithm which leads to $\mathbb{E}Q$. We will prove that, with such a choice of $m$, the complexity of the algorithm for $\mathbb{E}Q$ is significantly smaller than the complexity of the standard Monte Carlo method for $\mathbb{E}f(X_T^n)$.

Let us present the algorithm for $\mathbb{E}Q$.

3.2. *Central limit theorem.* In the following we assume that the parameters of the statistical Romberg method depend only on $n$. That is,

$$m = n^\beta, \qquad \beta \in (0,1), \qquad N_m = n^{\gamma_1}, \qquad \gamma_1 > 1, \qquad N_n = n^{\gamma_2}, \qquad \gamma_2 > 1.$$

We can now state the analogue of Theorem 3.1 in our setting. The statistical Romberg method approximates $\mathbb{E}f(X_T)$ by

$$V_n := \frac{1}{n^{\gamma_1}} \sum_{i=1}^{n^{\gamma_1}} f(\hat{X}_{T,i}^{n^\beta}) + \frac{1}{n^{\gamma_2}} \sum_{i=1}^{n^{\gamma_2}} f(X_{T,i}^n) - f(X_{T,i}^{n^\beta}),$$

where $\hat{X}_T^{n^\beta}$ is a second Euler scheme with step $T/n^\beta$ and such that the Brownian paths used for $X_T^n$ and $X_T^{n^\beta}$ have to be independent of the Brownian paths used in order to simulate $\hat{X}_T^{n^\beta}$. Here the quantity $\frac{1}{n^{\gamma_1}} \sum_{i=1}^{n^{\gamma_1}} f(\hat{X}_{T,i}^{n^\beta})$ must be viewed as an approximation for $M_m$.

THEOREM 3.2. *Let $f$ be an $\mathbb{R}^d$-valued function satisfying*

$$(15) \qquad |f(x) - f(y)| \le C(1 + |x|^p + |y|^p)|x - y| \qquad \textit{for some } C, p > 0.$$

*Assume that $\mathbb{P}(X_T \notin \mathcal{D}_f) = 0$, where $\mathcal{D}_f := \{x \in \mathbb{R}^d; f \text{ is differentiable at } x\}$, and that for some $\alpha \in [1/2, 1]$, we have*

$$(16) \qquad \lim_{n \to \infty} n^\alpha \varepsilon_n = C_f(T, \alpha).$$

*Then, for $\gamma_1 = 2\alpha$ and $\gamma_2 = 2\alpha - \beta$, we have*

$$n^\alpha(V_n - \mathbb{E}f(X_T)) \Longrightarrow \sigma_2 \tilde{G} + C_f(T, \alpha),$$

*with $\sigma_2^2 = \text{Var}(f(X_T)) + \tilde{\text{V}}\text{ar}(\nabla f(X_T)U_T)$ and $\tilde{G}$ a standard normal.*

LEMMA 3.1. *Under the assumptions of Theorem* 3.2, *for all* $\gamma > 0$,

$$\frac{1}{n^{(\gamma-\beta)/2}} \sum_{i=1}^{n^{\gamma}} f(X_{T,i}^{n^{\beta}}) - f(X_{T,i}) - \mathbb{E}(f(X_T^{n^{\beta}}) - f(X_T)) \Longrightarrow \mathcal{N}(0, \sigma_1^2), \tag{15}$$

*where*

$$\sigma_1^2 = \tilde{\mathrm{V}}\mathrm{ar}(\nabla f(X_T) \cdot U_T) \tag{16}$$

*and* $U$ *the process on* $\tilde{\mathcal{B}}$ *given by* (3).

PROOF. If we set

$$Z_T^{n^{\beta}} = f(X_T^{n^{\beta}}) - f(X_T) - \mathbb{E}(f(X_T^{n^{\beta}}) - f(X_T)),$$

then we have

$$\mathbb{E}\left[\exp\left(\frac{iu}{n^{((\gamma-\beta)/2)}} \sum_{k=1}^{n^{\gamma}} Z_{T,k}^{n^{\beta}}\right)\right] = \left[1 + \frac{1}{n^{\gamma}}\left(\frac{-u^2}{2} n^{\beta} \mathbb{E}\,|Z_T^{n^{\beta}}|^2 + n^{\gamma} \mathbb{E} C_n(\omega)\right)\right]^{n^{\gamma}},$$

*where*

$$|\mathbb{E} C_n(\omega)| \leq \frac{u^3}{6 n^{(3/2)(\gamma-\beta)}} \mathbb{E}\,|Z_T^{n^{\beta}}|^3.$$

Property (2) ensures the existence of a constant $K_3 > 0$ such that

$$|\mathbb{E} C_n(\omega)| \leq \frac{K_3 u^3}{6 n^{3\gamma/2}}.$$

We have, with probability 1,

$$n^{\beta/2}(f(X_T^{n^{\beta}}) - f(X_T)) = n^{\beta/2} \nabla f(X_T) \cdot U_T^{n^{\beta}} + R_n,$$

with

$$R_n = n^{\beta/2} |U_T^{n^{\beta}}| \varepsilon(X_T, U_T^{n^{\beta}}) \quad \text{and} \quad \varepsilon(X_T, U_T^{n^{\beta}}) \xrightarrow{\mathbb{P}} 0.$$

From the tightness of $(n^{\beta/2}|U_T^{n^{\beta}}|)_n$, it follows that $R_n \xrightarrow{\mathbb{P}} 0$, then, according to Lemma 2.1 and to Theorem 2.1,

$$n^{\beta/2}(f(X_T^{n^{\beta}}) - f(X_T)) \Longrightarrow \nabla f(X_T) \cdot U_T. \tag{17}$$

Using (13), it follows from property (2) that

$$\forall \varepsilon > 0 \qquad \sup_n \mathbb{E}|n^{\beta/2}(f(X_T^{n^{\beta}}) - f(X_T))|^{2+\varepsilon} < \infty.$$

Since $\mathbb{P}(X_T \notin \mathcal{D}_f) = 0$, we deduce, using (17), that

$$\mathbb{E}(n^{\beta/2}(f(X_T^{n^{\beta}}) - f(X_T)))^k \to \tilde{\mathbb{E}}(\nabla f(X_T) \cdot U_T)^k < \infty \qquad \text{with } k \in \{1, 2\}.$$

Consequently,

$$(18) \qquad n^{\beta}\mathbb{E}|Z_T^{n^{\beta}}|^2 \longrightarrow \tilde{\mathbb{V}}\mathrm{ar}(\nabla f(X_T) \cdot U_T) < \infty.$$

Since $\gamma > 0$, we see that $n^{\gamma}\mathbb{E}C_n(\omega) \to 0$ and we conclude that

$$\mathbb{E}\left[\exp\left(\frac{iu}{n^{((\gamma-\beta)/2)}}\sum_{k=1}^{n^{\gamma}} Z_{T,k}^{n^{\beta}}\right)\right] \longrightarrow \exp\left[\frac{-u^2}{2}\tilde{\mathbb{V}}\mathrm{ar}(\nabla f(X_T) \cdot U_T)\right],$$

which completes the proof.  $\square$

LEMMA 3.2.  *Under the assumptions of Theorem* 3.2, *for all* $\gamma > 0$,

$$\frac{1}{n^{((\gamma-\beta)/2)}}\sum_{k=1}^{n^{\gamma}} f(X_{T,k}^n) - f(X_{T,k}^{n^{\beta}}) - \mathbb{E}(f(X_T^n) - f(X_T^{n^{\beta}})) \Longrightarrow \mathcal{N}(0, \sigma_1^2),$$

*where* $\sigma_1^2 = \tilde{\mathbb{V}}\mathrm{ar}(\nabla f(X_T) \cdot U_T).$

PROOF.   We have

$$\frac{1}{n^{((\gamma-\beta)/2)}}\sum_{k=1}^{n^{\gamma}} f(X_{T,k}^n) - f(X_{T,k}^{n^{\beta}}) - \mathbb{E}(f(X_T^n) - f(X_T^{n^{\beta}}))$$

$$= \frac{1}{n^{((\gamma-\beta)/2)}}\sum_{k=1}^{n^{\gamma}} Z_{T,k}^n - \frac{1}{n^{((\gamma-\beta)/2)}}\sum_{k=1}^{n^{\gamma}} Z_{T,k}^{n^{\beta}}.$$

By (18), it follows that

$$(19) \qquad \mathbb{E}\left[\frac{1}{n^{((\gamma-\beta)/2)}}\sum_{k=1}^{n^{\gamma}} Z_{T,k}^n\right]^2 = n^{\beta}\mathbb{E}\left[Z_T^n\right]^2 \to 0.$$

The announced result follows from the above lemma.  $\square$

PROOF OF THEOREM 3.2.   For $\gamma_1 = 2\alpha, \gamma_2 = 2\alpha - \beta$, we have

$$n^{\alpha}(V_n - \mathbb{E}f(X_T)) = V_n^1 + V_n^2 + V_n^3,$$

where

$$(20) \qquad V_n^1 = \frac{1}{n^{\alpha}}\sum_{i=1}^{n^{2\alpha}} f(\hat{X}_{T,i}^{n^{\beta}}) - \mathbb{E}f(\hat{X}_T^{n^{\beta}}),$$

$$(21) \qquad V_n^2 = \frac{1}{n^{\alpha-\beta}}\sum_{i=1}^{n^{2\alpha-\beta}} f(X_{T,i}^n) - f(X_{T,i}^{n^{\beta}}) - \mathbb{E}(f(X_T^n) - f(X_T^{n^{\beta}})),$$

$$(22) \qquad V_n^3 = n^{\alpha}(\mathbb{E}f(X_T^n) - \mathbb{E}f(X_T)).$$

Properties (2) and (3) guarantee that the Lindeberg–Feller theorem applies here (same argument as in [5]). That is,

$$V_n^1 \Longrightarrow \mathcal{N}(0, \mathrm{Var}(f(X_T))).$$

On account of Lemma 3.2, it is obvious that

$$V_n^2 \Longrightarrow \mathcal{N}(0, \tilde{\mathrm{V}}\mathrm{ar}(\nabla f(X_T) \cdot U_T)).$$

Finally, by using the assumption (14), we complete the proof. □

3.3. *Complexity analysis.* As in the Monte Carlo case, we can interpret Theorem 3.2 as follows. For a total error of order $1/n^\alpha$, the minimal computational effort necessary to run the statistical Romberg algorithm applied to the Euler scheme (with step numbers $n$ and $m = n^\beta$) is obtained for

$$(23) \qquad N_m = n^2\alpha \quad \text{and} \quad N_n = n^{2\alpha - \beta}.$$

Since the only constraint on $\beta$ is that $\beta \in (0, 1)$, we will choose the optimal $\beta^\star$ minimizing the complexity of the statistical Romberg algorithm. The time complexity in the statistical Romberg method is given by

$$C_{\mathrm{SR}} = C \times (mN_m + (n + m)N_n) \qquad \text{with } C > 0$$
$$= C \times (n^{\beta + 2\alpha} + (n + n^\beta)n^{2\alpha - \beta}).$$

Simple calculations show that $\beta^\star = 1/2$ is the optimal choice which minimizes the time complexity.

So the optimal parameters in this case are the following:

$$m = n^{1/2}, \qquad N_m = n^{2\alpha} \quad \text{and} \quad N_n = n^{2\alpha - 1/2},$$

and the optimal complexity of the statistical Romberg method is given by

$$C_{\mathrm{SR}} \simeq C \times n^{2\alpha + 1/2}.$$

However, for the same error of order $1/n^\alpha$, we have shown that the optimal complexity of a Monte Carlo method was given by

$$C_{\mathrm{MC}} = C \times n^{2\alpha + 1},$$

which is clearly larger than $C_{\mathrm{SR}}$. So we deduce that the statistical Romberg method is more efficient.

**4. Statistical Romberg method and Asian options.** The payoff of an Asian option is related to the integral of the asset price process. Computing the price of an Asian option requires the discretization of the integral. The purpose of this section is to apply statistical Romberg extrapolation to the approximation of the integral and to carry on a complexity analysis in this context. This will lead us to prove a central limit theorem for the discretization error, which can be viewed as the analogue of Theorem 2.1 (see Theorem 4.1).

4.1. *Trapezoidal scheme.* Let $S$ be the process on the stochastic basis $\mathcal{B} = (\Omega, \mathcal{F}, (\mathcal{F}_t)_{t \geq 0}, P)$ satisfying

$$(24) \qquad \frac{dS_t}{S_t} = r\,dt + \sigma\,dW_t \qquad \text{with } t \in [0, T], T > 0,$$

where $\sigma$ and $r$ are real constants, with $\sigma > 0$ and $(W_t)_{t \in [0,T]}$ is a standard Brownian motion on $\mathcal{B}$. The solution of the last equation is given by

$$(25) \qquad S_t = S_0 \exp\left(\left(r - \frac{\sigma^2}{2}\right)t + \sigma W_t\right).$$

We set

$$I_T = \frac{1}{T} \int_0^T S_u\,du.$$

Let $f$ be a given real valued function. Our aim will be to evaluate

$$\Pi(S, T) = e^{-rT}\mathbb{E}f(S_T, I_T).$$

In a financial setting, if $f(x, y) = (y - K)_+$, $\Pi(S, T)$ is the price of an Asian call option with fixed strike $K$. In this case there is no explicit formula that gives the real price. So, the computation of this price, by a probabilistic method, requires a discretization of the integral $I_T$. There are several approximation schemes used in practice. One of the most efficient is the trapezoidal scheme defined by

$$(26) \qquad I_T^n = \frac{\delta}{T} \sum_{k=1}^n S_{t_{k-1}} \left(1 + \frac{r\delta}{2} + \sigma \frac{W_{t_k} - W_{t_{k-1}}}{2}\right),$$

where $\delta = \frac{T}{n}$ and $t_k = \frac{Tk}{n} = \delta k$. We call it trapezoidal because it is closely related to the trapezoidal approximation of the integral

$$\mathbb{E}\left(I_T^n - \frac{\delta}{T} \sum_{k=1}^n \frac{S_{t_{k-1}} + S_{t_k}}{2}\right)^2 = O\left(\frac{1}{n^3}\right).$$

Note that $S_{t_k}$ has an explicit expression so we can simulate it without discretizing the SDE. The following result is proved in [17].

PROPOSITION 4.1. *With the above notation, there exists a nondecreasing map $K(T)$ such that, $\forall p > 0$,*

$$\left[\mathbb{E}\left(\sup_{t \in [0,T]} |I_t^n - I_t|^{2p}\right)\right]^{1/(2p)} \leq \frac{K(T)}{n}.$$

4.2. *Stable convergence of the trapezoidal scheme error.*  In the following we prove a functional CLT analogous to Jacod and Protter's theorem (see Theorem 2.1 above).

THEOREM 4.1.  *Let* $J_t = \frac{1}{T} \int_0^t S_u \, du, t \in [0, T]$, *and* $J_n$ *be the trapezoidal discretization associated with* $J$:

$$J_t^n := \frac{\delta}{T} \sum_{k=1}^{[t/\delta]} S_{t_{k-1}} \left( 1 + \frac{r\delta}{2} + \sigma \frac{W_{t_k} - W_{t_{k-1}}}{2} \right).$$

*We have*

(27)                              $$n(J - J^n) \stackrel{stably}{\Longrightarrow} \chi,$$

*where* $\chi$ *is the process defined by*

$$\chi_t = \frac{\sigma}{2\sqrt{3}} \int_0^t S_s \, dB'_s,$$

*where* $B'$ *is a standard Brownian motion on an extension* $\hat{\mathcal{B}}$ *of* $\mathcal{B}$, *which is independent of* $W$.

We have the following elementary lemma.

LEMMA 4.1.  *If* $H$ *is deterministic satisfying* $\int_0^t H_s^2 \, ds < \infty$, *then we have*

$$\int_0^t \frac{\tau_\delta(s)}{\delta} H_s \, ds \longrightarrow \frac{1}{2} \int_0^t H_s \, ds$$

*and*

$$\int_0^t \left( 2\frac{\tau_\delta(s)}{\delta} - 1 \right)^2 H_s \, ds \longrightarrow \frac{1}{3} \int_0^t H_s \, ds,$$

*with* $\tau_\delta(s) = t \wedge ([s/\delta]\delta + \delta) - s$.

PROOF.   We sketch the proof for completeness. By a density argument, we may assume that $H$ is piecewise constant: $H_s = c_i$ for $T_{i-1} < s < T_i$, where $0 = T_0 < \cdots < T_k = T$ and $(c_i)$ are constants for $i = 0, \ldots, k$. It follows that

$$\int_0^t \frac{\tau_\delta(s)}{\delta} H_s \, ds = \sum_{i=1}^k \frac{c_i}{\delta} \int_{T_i}^{T_{i+1}} \tau_\delta(s) \, ds \to \sum_{i=1}^k \frac{c_i}{2} (T_{i+1} - T_i) = \int_0^t \frac{H_s}{2} \, ds,$$

since it is easy to check that

$$\frac{1}{\delta} \int_x^y \tau_\delta(s) \, ds \longrightarrow \frac{y - x}{2} \qquad \text{as } n \longrightarrow \infty.$$

The second assertion is obtained in the same way, but by using that

$$\frac{1}{\delta^2}\int_x^y\left(2\frac{\tau_\delta(s)}{\delta}-\delta\right)^2 ds \longrightarrow \frac{y-x}{3} \qquad \text{as } n \longrightarrow \infty,$$

which completes the proof.  □

PROOF OF THEOREM 4.1.   We have

$$
n(J_t - J_t^n) = \frac{n}{T}\int_0^t S_u\,du - n\Bigg(\frac{\delta}{T}\sum_{k=1}^{[t/\delta]}S_{t_{k-1}} + \frac{\delta^2 r}{2T}\sum_{k=1}^{[t/\delta]}S_{t_{k-1}}
$$

$$
(28) \qquad\qquad\qquad\qquad + \frac{\delta\sigma}{2T}\sum_{k=1}^{[t/\delta]}S_{t_{k-1}}(W_{t_k}-W_{t_{k-1}})\Bigg)
$$

$$
+ \frac{n\delta}{T}S_{[t/\delta]\delta}\Big(1 + \frac{r}{2}(t-[t/\delta]\delta) + \frac{\sigma}{2}(W_t - W_{[t/\delta]\delta})\Big).
$$

It follows that

$$
(29) \qquad n(J_t - J_t^n) = A_t^\delta - \frac{r}{2}\int_0^t S_{[u/\delta]\delta}\,du - \frac{\sigma}{2}\int_0^t S_{[u/\delta]\delta}\,dW_u,
$$

with

$$
A_t^\delta = \frac{1}{\delta}\int_0^t (S_u - S_{[u/\delta]\delta})\,du.
$$

Note that, by using (24), we obtain

$$
A_t^\delta = A_t^{\delta,1} + A_t^{\delta,2},
$$

with

$$
A_t^{\delta,1} = \frac{r}{\delta}\int_0^t\int_{[u/\delta]\delta}^u S_s\,ds\,du
$$

and

$$
A_t^{\delta,2} = \frac{\sigma}{\delta}\int_0^t\int_{[u/\delta]\delta}^u S_s\,dW_s\,du,
$$

and we have

$$
A_t^{\delta,1} = \frac{r}{\delta}\sum_{k=1}^{[t/\delta]}\int_{t_{k-1}}^{t_k}\int_{t_{k-1}}^u S_s\,ds\,du + \frac{r}{\delta}\int_{[t/\delta]\delta}^t\int_{[t/\delta]\delta}^u S_s\,ds\,du,
$$

$$
= \frac{r}{\delta}\int_0^t (t\wedge([s/\delta]\delta+\delta)-s)S_s\,ds,
$$

$$
= \frac{r}{\delta}\int_0^t \tau_\delta(s)S_s\,ds.
$$

In the same manner we can see that

$$A_t^{\delta,2} = \frac{\sigma}{\delta} \sum_{k=1}^{[t/\delta]} \int_{t_{k-1}}^{t_k} \int_{t_{k-1}}^u S_s \, dW_s \, du + \frac{\sigma}{\delta} \int_{[t/\delta]\delta}^t \int_{[t/\delta]\delta}^u S_s \, dW_s \, du.$$

The integration by parts formula gives

$$\int_{t_{k-1}}^{t_k} \int_{t_{k-1}}^u S_s \, dW_s \, du = (t_k - t_{k-1}) \int_{t_{k-1}}^{t_k} S_s \, dW_s - \int_{t_{k-1}}^{t_k} (s - t_k) S_s \, dW_s$$

$$= \int_{t_{k-1}}^{t_k} ([s/\delta]\delta + \delta - s) S_s \, dW_s$$

and

$$\int_{[t/\delta]\delta}^t \int_{[t/\delta]\delta}^u S_s \, dW_s \, du = (t - [t/\delta]\delta) \int_{[t/\delta]\delta}^t S_s \, dW_s$$

$$- \int_{[t/\delta]\delta}^t (s - [t/\delta]\delta) S_s \, dW_s$$

$$= \int_{[t/\delta]\delta}^t (t - s) S_s \, dW_s.$$

We deduce that

$$A_t^{\delta,2} = \frac{\sigma}{\delta} \sum_{k=1}^{[t/\delta]} \int_{t_{k-1}}^{t_k} ([s/\delta]\delta + \delta - s) S_s \, dW_s + \frac{\sigma}{\delta} \int_{[t/\delta]\delta}^t (t - s) S_s \, dW_s,$$

$$= \frac{\sigma}{\delta} \int_0^t (t \wedge ([s/\delta]\delta + \delta) - s) S_s \, dW_s,$$

$$= \frac{\sigma}{\delta} \int_0^t \tau_\delta(s) S_s \, dW_s.$$

It follows that

$$n(J_t - J_t^n) = \frac{r}{\delta} \int_0^t \tau_\delta(s) S_s \, ds + \frac{\sigma}{\delta} \int_0^t \tau_\delta(s) S_s \, dW_s$$

$$- \frac{r}{2} \int_0^t S_{[s/\delta]\delta} \, ds - \frac{\sigma}{2} \int_0^t S_{[s/\delta]\delta} \, dW_s.$$

We deduce that

$$n(J_t - J_t^n) = B_t^\delta + \chi_t^\delta + C_t^\delta,$$

with

$$B_t^\delta = \frac{r}{2} \int_0^t \left( 2 \frac{\tau_\delta(s)}{\delta} - 1 \right) S_s \, ds,$$

$$\chi_t^\delta = \frac{\sigma}{2} \int_0^t \left( 2\frac{\tau_\delta(s)}{\delta} - 1 \right) S_s \, dW_s,$$

$$C_t^\delta = \frac{r}{2} \int_0^t (S_s - S_{[s/\delta]\delta}) \, ds + \frac{\sigma}{2} \int_0^t (S_s - S_{[s/\delta]\delta}) \, dW_s.$$

According to the above lemma, we obtain that $\sup_{t \in [0,T]} B_t^\delta \to 0$ a.s. It is obvious that $\sup_{t \in [0,T]} C_t^\delta \xrightarrow{L^2} 0$. The only point remaining concerns the behavior of $\chi^\delta$.

In virtue of Theorem 2.1 of [7], if we prove that, for all $t \in [0, T]$, we have

$$\langle \chi^\delta, \chi^\delta \rangle_t \xrightarrow{\mathbb{P}} \frac{\sigma^2}{12} \int_0^t S_s^2 \, ds, \qquad \langle \chi^\delta, W \rangle_t \xrightarrow{\mathbb{P}} 0,$$

then the process $\chi^\delta$ will converge stably in law to the process $\chi$ of (27). According to the above lemma, we have

$$\langle \chi^\delta, \chi^\delta \rangle_t = \frac{\sigma^2}{4} \int_0^t \left( 2\frac{\tau_\delta(s)}{\delta} - 1 \right)^2 S_s^2 \, dW_s \longrightarrow \frac{\sigma^2}{12} \int_0^t S_s^2 \, ds \qquad \text{a.s.}$$

and

$$\langle \chi^\delta, W \rangle_t = \frac{\sigma}{2} \int_0^t \left( 2\frac{\tau_\delta(s)}{\delta} - 1 \right) S_s \, ds \longrightarrow 0 \qquad \text{a.s.,}$$

which completes the proof. $\quad\square$

4.3. *Statistical Romberg method and CLT.* In order to evaluate $e^{-rT}\mathbb{E}f(S_T, I_T)$, we use the idea of statistical Romberg approximation:

- compute an approximation $E_n^1$ of $e^{-rT}\mathbb{E}f(S_T, \hat{I}_T^{n^\beta})$ by a Monte Carlo method

$$E_n^1 = \frac{e^{-rT}}{n^{\gamma_1}} \sum_{i=1}^{n^{\gamma_1}} f(\hat{S}_{T,i}, \hat{I}_{T,i}^{n^\beta}),$$

- compute

$$E_n^2 = \frac{e^{-rT}}{n^{\gamma_2}} \sum_{i=1}^{n^{\gamma_2}} f(S_{T,i}, I_{T,i}^n) - f(S_{T,i}, I_{T,i}^{n^\beta}).$$

Recall that the samples used, in order to construct

$$(\hat{S}_{T,i}, \hat{I}_{T,i}^{n^\beta}) \quad \text{and} \quad ((S_{T,i}, I_{T,i}^n), (S_{T,i}, I_{T,i}^{n^\beta})),$$

are independent. The question now is how to choose $\beta, \gamma_1$ and $\gamma_2$. Admittedly, we can choose the optimal parameters given in the Euler scheme case, but the following result proves that, in the specific case of the trapezoidal scheme, the optimal parameters are different.

THEOREM 4.2. *Let $f$ be an $\mathbb{R}^2$-valued function satisfying*

$$(30) \qquad |f(x, y_1) - f(x, y_2)| \leq C(1 + |x|^p + |y_1|^p + |y_2|^p)|y_1 - y_2|$$

$$\text{for some } C, p > 0.$$

*Assume that $\mathbb{P}((S_T, I_T) \notin \mathscr{D}_f) = 0$, with*

$$\mathscr{D}_f := \{(x, y) \in \mathbb{R}^+ \times \mathbb{R}^+ \,|\, \partial_2 f(x, y) \text{ exists}\},$$

*where $\partial_2 f$ denotes the partial derivative of $f$ w.r.t. $y$. Then for all $\beta \in (0, 1)$, if $\gamma_1 = 2$ and $\gamma_2 = 2 - 2\beta$, we have*

$$n(E_n - \mathbb{E}f(S_T, I_T)) \Longrightarrow \tilde{\sigma}_2 \hat{G} + \hat{\mathbb{E}}(\partial_2 f(S_T, I_T)\chi_T),$$

*where $\hat{\sigma}_2^2 = \mathrm{Var}(f(S_T, I_T)) + \hat{\mathrm{V}}\mathrm{ar}(\partial_2 f(S_T, I_T)\chi_T)$, $\chi$ is the limit process on $\hat{\mathcal{B}}$ given in Theorem 4.1, and $\hat{G}$ a standard normal.*

REMARK 4.1. The assumptions on $f$ in the above theorem are satisfied in the case of typical Asian options:

$$f(x, y) = (y - K)_+, \qquad f(x, y) = (K - y)_+, \qquad f(x, y) = (y - x)_+.$$

LEMMA 4.2. *Under the assumptions of Theorem 4.2, for all $\gamma > 0$,*

$$\frac{1}{n^{((\gamma/2) - \beta)}} \sum_{k=1}^{n^\gamma} f(S_{T,k}, I_{T,k}) - f(S_{T,k}, I_{T,k}^{n^\beta})$$

$$- \mathbb{E}(f(S_T, I_T) - f(S_T, I_T^{n^\beta})) \Longrightarrow \mathcal{N}(0, \hat{\sigma}_1^2),$$

*where $\hat{\sigma}_1^2 = \hat{\mathrm{V}}\mathrm{ar}(\partial_2 f(S_T, I_T)\chi_T)$.*

PROOF. If we set

$$H_T^{n^\beta} := f(S_T, I_T) - f(S_T, I_T^{n^\beta}) - \mathbb{E}(f(S_T, I_T) - f(S_T, I_T^{n^\beta})),$$

then we have

$$\mathbb{E}\left[\exp\left(\frac{iu}{n^{((\gamma/2) - \beta)}} \sum_{k=1}^{n^\gamma} H_{T,k}^{n^\beta}\right)\right]$$

$$(30) \qquad = \left[1 + \frac{1}{n^\gamma}\left(\frac{-u^2}{2} \mathrm{Var}[n^\beta(f(S_T, I_T) - f(S_T, I_T^{n^\beta}))]\right.\right.$$

$$\left.\left. + n^\gamma \mathbb{E}C_n'(\omega)\right)\right]^{n^\gamma},$$

with

$$|EC_n'(\omega)| \leq \frac{u^3}{6n^{3((\gamma/2) - \beta)}} \mathbb{E}|H_T^{n^\beta}|^3.$$

Proposition 4.1 ensures the existence of a constant $K(T) > 0$ such that

$$|EC'_n(\omega)| \leq \frac{K(T)u^3}{6n^{3\gamma/2}}.$$

We have, with probability 1,

$$n^\beta(f(S_T, I_T^{n^\beta}) - f(S_T, I_T)) = n^\beta \partial_2 f(S_T, I_T)(I_T^{n^\beta} - I_T) + R_n,$$

with

$$R_n = n^\beta |I_T^{n^\beta} - I_T| \varepsilon(S_T, I_T, I_T^{n^\beta}) \quad \text{and} \quad \varepsilon(S_T, I_T, I_T^{n^\beta}) \xrightarrow{\mathbb{P}} 0.$$

From the tightness of $(n^\beta |I_T^{n^\beta} - I_T|)_n$, it follows that $R_n \xrightarrow{\mathbb{P}} 0$. Consequently, according to Lemma 2.1 and to Theorem 4.1, we obtain that

$$(31) \qquad n^\beta(f(S_T, I_T^{n^\beta}) - f(S_T, I_T)) \Longrightarrow \partial_2 f(S_T, I_T)\chi_T.$$

With our assumption (30) on $f$, it follows from Proposition 4.1 that

$$\sup_n \mathbb{E}|n^\beta(f(S_T, I_T) - f(S_T, I_T^{n^\beta}))|^{2+\varepsilon} < \infty \qquad \text{with } \varepsilon > 0,$$

so we obtain

$$(32) \qquad \begin{aligned} &\mathbb{E}(n^\beta(f(S_T, I_T) - f(S_T, I_T^{n^\beta})))^k \\ &\qquad \to \hat{\mathbb{E}}(\partial_2 f(S_T, I_T)\chi_T)^k < \infty \qquad \forall 0 < k \leq 2. \end{aligned}$$

Hence, we deduce that

$$\mathbb{E}\left[\exp\left(\frac{iu}{n^{((\gamma/2)-\beta)}} \sum_{k=1}^{n^\gamma} H_{T,k}^{n^\beta}\right)\right] \longrightarrow \exp\left[\frac{-u^2}{2}\hat{V}\mathrm{ar}(\partial_2 f(S_T, I_T)\chi_T)\right],$$

which completes the proof.  □

REMARK 4.2. It follows from the proof of the above lemma that

$$(33) \qquad \lim_{n \to \infty} n\mathbb{E}(f(S_T, I_T^n) - f(S_T, I_T)) = \widetilde{\mathbb{E}(\partial_2 f(S_T, I_T)\chi_T)}.$$

This gives us an expansion of the discretization error in our setting. Note that similar expansions are given in  [17] for less regular functions. The advantage of our approach is that we do not need the associate PDE, so that our expansion is more explicit.

The proof of the following result is a consequence of the above lemma.

LEMMA 4.3. *Under the assumptions of Theorem 4.2 and for all $\gamma > 0$, we have*

$$\frac{1}{n^{((\gamma/2)-\beta)}} \sum_{k=1}^{n^{\gamma}} f(S_{T,k}, I_{T,k}^n) - f(S_{T,k}, I_{T,k}^{n^{\beta}})$$

$$- \mathbb{E}(f(S_T, I_T^n) - f(S_T, I_T^{n^{\beta}})) \Longrightarrow \mathcal{N}(0, \hat{\sigma}_1^2),$$

*with $\hat{\sigma}_1^2 = \tilde{\mathrm{Var}}(\partial_2 f(S_T, I_T)\chi_T)$.*

Using Lemma 4.3, Theorem 4.2 can be proved in much the same way as Theorem 3.2. Equality (33) will be used instead of the assumption (14) on $\lim_{n\to\infty} n^{\alpha}\varepsilon_n$ in Theorem 3.2.

4.4. *Complexity analysis.* As in the Euler scheme case, one can interpret the above theorem in the following way. For a total error of order $1/n$, the minimal computational effort necessary to run the statistical Romberg method applied to the trapezoidal scheme (with step numbers $n$ and $m = n^{\beta}$) is obtained for

(34) $$N_m = n^2 \quad \text{and} \quad N_n = n^{2-2\beta}.$$

Since the only restriction on $\beta$ is that $\beta \in (0,1)$, we will choose the optimal $\beta^{\star}$ minimizing the complexity of the statistical Romberg algorithm. In this case the time complexity in the statistical Romberg method is given by

$$C_{\mathrm{SR}} = C \times (mN_m + (n+m)N_n) \qquad \text{with } C > 0,$$
$$= C \times (n^{\beta+2} + (n+n^{\beta})n^{2-2\beta}).$$

Simple calculations show that $\beta^{\star} \simeq 1/3$ is the optimal choice which minimizes the time complexity.

So the optimal parameters in this case are

$$m = n^{1/3}, \qquad N_m = 2 \quad \text{and} \quad N_n = n^{4/3},$$

and the optimal complexity of the statistical Romberg method is given by

$$C_{\mathrm{SR}} \simeq C \times n^{7/3}.$$

But according to Proposition 4.1, Theorem 3.1 remains valid if we change the Euler scheme by the trapezoidal one. Hence, for the same error of order $1/n$, the optimal complexity of a Monte Carlo method applied to the trapezoidal scheme with step number $n$ is given by

$$C_{\mathrm{MC}} = C \times n^3,$$

which is clearly larger than $C_{\mathrm{SR}}$. So we deduce that the statistical Romberg method is more efficient.

**5. Numerical tests and results.** We test the efficiency of the statistical Romberg method to reduce the time complexity for the degenerate two-dimensional diffusion given in the example of Section 2 (tests concerning Asian options are given in [10]).

Consider the bi-dimensional diffusion

$$Z_t = (\cos(\theta + W_t), \sin(\theta + W_t)), \qquad \theta \in [0, 2\pi],$$

and the map

$$f_\alpha : z = (x, y) \mapsto ||z|^2 - 1|^{2\alpha}, \qquad \alpha \in [1/2, 1],$$

given in the example of Section 2. Note that in this case

$$\text{Var}(f_\alpha(Z_T)) = 0,$$

and this is because the bi-dimensional diffusion $Z = (X, Y)$ lives on the unit circle. Consequently, in order to obtain the optimal parameters given by Theorems 3.1 and 3.2, we consider $g_\alpha(x, y) = f_\alpha(x, y) + x$, instead of $f_\alpha(x, y)$. This choice leads us to a nonvanishing variance of $g_\alpha(Z_T)$ and to a discretization error which is of order $1/n^\alpha$, $\alpha \in [1/2, 1]$. In the following we set

- *MC method*: the algorithm using a Monte Carlo method to approximate $\mathbb{E}(g_\alpha(Z_T))$ by

$$\frac{1}{N} \sum_{i=1}^{N} g_\alpha(Z_{T,i}^n).$$

- *SR method*: the algorithm using a statistical Romberg method to approximate $\mathbb{E}(g_\alpha(Z_T))$ by

$$(35) \qquad \frac{1}{N_n} \sum_{i=1}^{N_n} [g_\alpha(Z_{T,i}^n) - g_\alpha(Z_{T,i}^m)] + \frac{1}{N_m} \sum_{i=1}^{N_m} g_\alpha(\hat{Z}_{T,i}^m).$$

To compare both methods, we use the methodology proposed by Broadie and Detemple [4]. Their idea is that, for a given set of parameters of the concerned diffusion, one of both algorithms will give better results. So they propose to test the algorithm on a large set of parameters chosen randomly. Proceeding along this line, we produce randomly $M = 200$ values for $Z_0 = (X_0, Y_0)$. Then, for each method, we compute the speed and an error measure. Speed is measured by the number of simulated values computed per second (the computations were done on a PC with a 2.00 GHz Pentium 4 processor). The error measure is given by the root-mean-squared error, which is defined by

$$(36) \qquad RMS = \sqrt{\frac{1}{M} \sum_{i=1}^{M} (\text{Real value} - \text{Simulated value})^2},$$

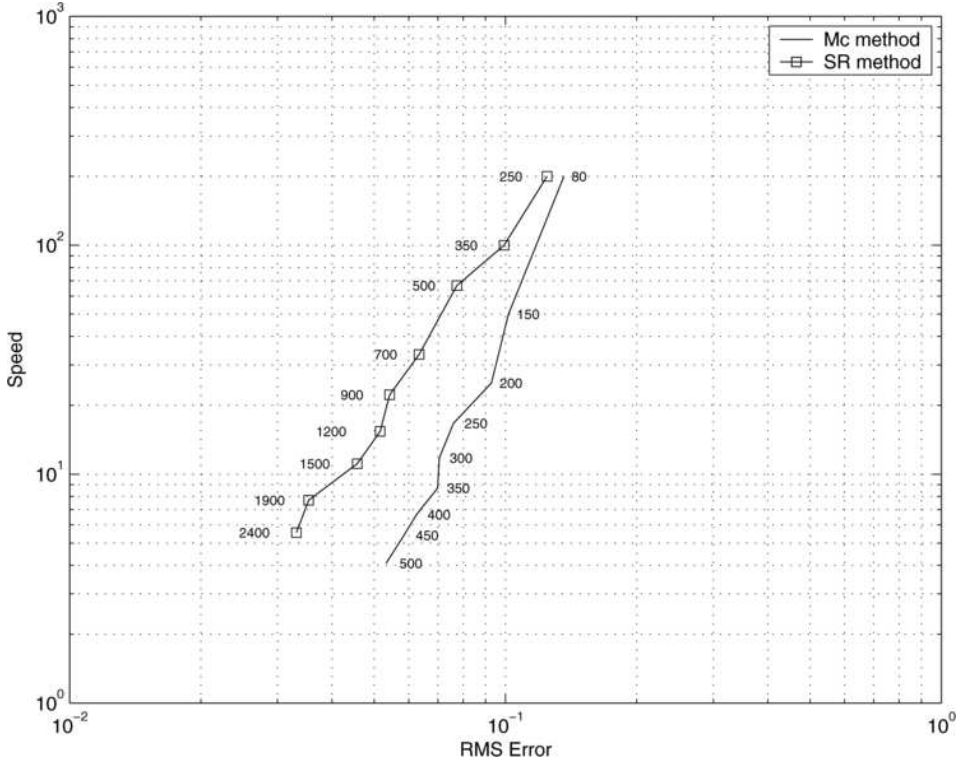

Fɪɢ. 1. *Speed versus RMS-error.*

and the real value is given by the formula $\mathbb{E}(g_\alpha(Z_T)) = \cos(\theta - T/2)$.

Our algorithm proceeds as follows. We fix the number of steps, say, $n = 80$, in the Euler scheme. We compute by the Monte Carlo method (resp., the statistical Romberg method) the 200 simulated values. Then, we produce according to (36) RMS$^{(MC,n=200)}$ and we compute Speed$^{(MC,n=200)}$, so we have a couple of points

$$(\text{RMS}^{(MC,n=200)}, \text{Speed}^{(MC,n=200)}).$$

This point is plotted on Figure 1. So, each given $n$ produces a couple of points: $(MC^{(MC,n)}, \text{Speed}^{(MC,n)})$. This gives the continuous curve in Figure 1. The line marked by squares is produced in the same way, using the statistical Romberg method.

Tests are done for $\alpha = 1/2$. Note that, although $g_\alpha$ is not $\mathscr{C}^1$ in that case, Theorem 3.2 can be extended to this specific example, as we can handle the difference $g_{1/2}(Z_T^{\sqrt{n}}) - g_{1/2}(Z_T)$.

Let us now interpret the curves. The fact that the SR curve is higher than the MC one means that, given an error $\varepsilon$, the number of values computed in one second (with this error) by the statistical Romberg algorithm is larger than the number of values computed by the Monte Carlo method. Note anyway that, for a large $\varepsilon$ (which corresponds to a small number of steps $n$), the differences between the two methods is less important. But as $\varepsilon$ becomes small ($n$ becomes large), the difference becomes more significant. In Table 1 we compare the speed of the Monte Carlo method and the speed of the statistical Romberg one, for a fixed RMS-error. We note that, by using the statistical Romberg method and for an RMS-error fixed at $10^{-1}$, one increases the speed by a factor of 2.26. For a small RMS-error fixed at $6 \cdot 10^{-2}$, the speed gain reaches a factor of 4.96.

TABLE 1
*Time complexity reduction for $\alpha = 1/2$*

| RMS relative error | MC method speed | SR method speed |
|---|---|---|
| $10^{-1}$ | 45.405 | 102.718 |
| $9 \cdot 10^{-2}$ | 23.599 | 85.876 |
| $8 \cdot 10^{-2}$ | 18.647 | 70.712 |
| $7 \cdot 10^{-2}$ | 9.219 | 49.179 |
| $6 \cdot 10^{-2}$ | 5.883 | 29.234 |

**6. Conclusion.** The statistical Romberg algorithm is a method that can be used in a general framework: as soon as we use a discretization scheme for the diffusion $(X_t)_{0 \leq t \leq T}$ in order to compute quantities such as $\mathbb{E}f(X_T)$, we can implement the statistical Romberg algorithm. And this is worth it because it is more efficient than a classic Monte Carlo method.

In financial applications, it is sometimes essential to be able to price a given product on the market as soon as possible and this by setting a margin of error that one can tolerate. In this case the statistical Romberg method equipped with its parameters allowing complexity reduction is faster than the standard Monte Carlo one.

## REFERENCES

[1] ALDOUS, D. J. and EAGLESON, G. K. (1978). On mixing and stability of limit theorems. *Ann. Probab.* **6** 325–331. MR517416

[2] BALLY, V. and TALAY, D. (1996). The law of the Euler scheme for stochastic differential equations. I. Convergence rate of the distribution function. *Probab. Theory Related Fields* **104** 43–60. MR1367666

[3] BOYLE, P., BROADIE, M. and GLASSERMAN, P. (1997). Monte Carlo methods for security pricing. *J. Econom. Dynam. Control* **21** 1267–1321. MR1470283

[4] BROADIE, M. and DETEMPLE, J. (1997). Recent advances in numerical methods for pricing derivative securities. In *Numerical Methods in Finance* 43–66. Cambridge Univ. Press. MR1470508

[5] DUFFIE, D. and GLYNN, P. (1995). Efficient Monte Carlo simulation of security prices. *Ann. Appl. Probab.* **5** 897–905. MR1384358

[6] Faure, O. (1992). Simulation du mouvement Brownien et des diffusions. Ph.D. thesis, Ecole Nationale des Ponts et Chaussées.

[7] Jacod, J. (1997). On continuous conditional Gaussian martingales and stable convergence in law. *Séminaire de Probabilités XXXI. Lecture Notes in Math.* **1655** 232–246. Springer, Berlin. MR1478732

[8] Jacod, J., Kurtz, T. G., Méléard, S. and Protter, P. (2005). The approximate Euler method for Lévy driven stochastic differential equations. *Ann. Inst. H. Poincaré. Special Issue Devoted to the Memory of P. A. Meyer* **41** 523–558. MR2139032

[9] Jacod, J. and Protter, P. (1998). Asymptotic error distributions for the Euler method for stochastic differential equations. *Ann. Probab.* **26** 267–307. MR1617049

[10] Kebaier, A. (2004). Romberg extrapolation, variance reduction and applications to option pricing. Preprint.

[11] Kurtz, T. G. and Protter, P. (1991). Wong–Zakai corrections, random evolutions, and simulation schemes for SDEs. In *Stochastic Analysis* 331–346. Academic Press, Boston, MA. MR1119837

[12] Kurtz, T. G. and Protter, P. (1999). Weak error estimates for simulation schemes for SDEs.

[13] Lamberton, D. and Lapeyre, B. (1997). *Introduction au calcul stochastique appliqué à la finance*, 2nd ed. Ellipses, Paris. MR1607509

[14] Protter, P. (1990). *Stochastic Integration and Differential Equations.* Springer, Berlin. MR1037262

[15] Rényi, A. (1963). On stable sequences of events. *Sankhyā Ser. A* **25** 293–302. MR170385

[16] Talay, D. and Tubaro, L. (1990). Expansion of the global error for numerical schemes solving stochastic differential equations. *Stochastic Anal. Appl.* **8** 483–509. MR1091544

[17] Temam, E. (2001). Couverture approché d'options exotiques. Pricing des options Asiatiques. Ph.D. thesis, Univ. Paris VI.

Laboratoire d'Analyse et de Mathématiques
UMR 8050
Université de Marne-La-Vallée
Cité Descartes, 5 Boulevard Descartes
Champs-sur-Marne
F-77454 Marne-La-Vallée Cedex 2
France
e-mail: ahmed.kbaier@univ-mlv.fr